\documentclass[a4paper,12pt]{article}

\usepackage{amsmath,amssymb,amsthm}
\usepackage{bbm, mathtools}
\usepackage{color}

\usepackage{cite}           
\newtheorem{theorem}{Theorem}
\newtheorem{lemma}{Lemma}
\newtheorem{corollary}{Corollary}
\newtheorem{proposition}[theorem]{Proposition}

\newcommand{\Z}{{\mathbb{Z}}}

\newcommand{\R}{{\mathbb{R}}}

\newcommand{\m}{{\tilde{m}}}
\def\={\stackrel {\rm def}  {=}}
\def\di={\overset{\text{${\mathcal{D} }$}} =}

\def\e{\operatorname{e}}

\title{Optimal Littlewood-Offord inequalities in groups}
\author{T. Ju\v skevi\v cius$^{1}$, G. \v Semetulskis$^{2}$}
\date{}
\begin{document}
\maketitle

\footnotetext[1] {Vilnius University, Institute of Mathematics and Informatics, Vilnius, Lithuania, email - tomas.juskevicius@gmail.com.}
\footnotetext[2] {University of Vilnius, Vilnius, Lithuania, email - grazvydas.semetulskis@gmail.com.}

\begin{abstract}
We prove several Littlewood-Offord type inequalities for arbitrary groups. In
groups having elements of finite order the worst case scenario is provided by
the simple random walk on a cyclic subgroup. The inequalities we obtain are
optimal if the underlying group contains an element of a certain order. 
It turns out that for torsion-free groups Erd\H{o}s's bound still holds. Our results strengthen and generalize some very recent results by Tiep and Vu.
\end{abstract}

\section{Introduction}
Let $V_n=\lbrace g_1,\ldots, g_n\rbrace$ be a multiset of non-identify elements of an arbitrary group $G$. Consider a collection of independent random variables $X_i$ that are each distributed on a two point set $\lbrace g^{-1}_i, g_{i} \rbrace$ and define the quantity
$$\rho(V_n)=  \sup_{g\in G}\mathbb{P}(X_1\ast\cdots\ast X_n=g).$$ In the case $G=\mathbb{R}$ the latter quantity is the maximum probability of the sum $X_1+\ldots+X_n$. 
Whenever $G=\R$, $\Z_m$ and $g_i = 1$ we shall adopt the convention to write 
$\varepsilon_i$ instead of the random variable $X_i$.

\quad Investigating random polynomials Littlewood and Offord \cite{LO} proved an 
almost optimal bound for the probability that a sum of random signs with
non-zero weights hits a point. To be more precise, using harmonic analysis they proved that in the case $G=\mathbb{R}$ we have

$$\rho (V_n)= O(n^{-1/2}\log n).$$

Erd\H{o}s \cite{ELO}, using Sperner's theorem from finite set combinatorics, showed that, actually,
$$\rho (V_n)\leq \frac{\binom{n}{\left\lfloor n/2\right\rfloor}}{2^n}.$$

\quad This bound is optimal as can be seen by taking $g_i=1$ in $V_n$. 
In this case we have 
$$
\rho(V_n) = \mathbb{P}(\varepsilon_1+\cdots+ \varepsilon_n \in \{0,1\}) = \frac{\binom{n}{\left\lfloor n/2\right\rfloor}}{2^n}.
$$
Answering a question of Erd\H{o}s, Kleitman \cite{KLO} used an ingenious
induction to show that the latter bound still holds for $g_i$ lying in
an arbitrary normed space. See also \cite{bollobas} for a very nice exposition of
Kleitman's beautiful argument. 
Griggs used a similar approach in \cite{Griggs} as in Erd\H{o}s's seminal paper 
\cite{ELO} to obtain the best possible result in $\Z_m$.

 More recently Tiep and Vu \cite{Tiep} investigated the same question for certain matrix groups and obtained results that are sharp up to a constant factor. 
To be more precise, let $m,k,n \geq 2$ be integers and $G=GL_{k}(\mathbb{C})$. 
Let $V_n=\lbrace g_1,\ldots, g_n \rbrace$ be a multiset of elements in $G$, each of which has order at least $m$. In this case they have obtained the bound

\begin{equation}\label{TiepVU}
\rho(V_n)\leq 141\max\lbrace\frac{1}{m},\frac{1}{\sqrt{n}}\rbrace.
\end{equation}
Furthermore, they have also established the same bound for $GL_{k}(p)$. 

Let us explain the meaning of the two terms in the upper bound given in (\ref{TiepVU}).
Take some element $g$ in $G$ of order $m$ and consider the multiset $V_n=\lbrace g,\ldots, g \rbrace$. 
Let us for the simplicity assume that $m$ is odd. In this setup the random
variable $S_k = X_1 \ast \dots \ast X_k$ is just the simple random walk on a subgroup of
$G$ that is isomorphic to $\Z_m$. It is a well
known fact that the distribution of $S_n$ is asymptotically uniform,
which accounts for the $\frac{1}{m}$ term in ($\ref{TiepVU}$).
For $n < m$ the point masses of $S_n$ are just the usual binomial probabilities
$\binom{n}{\left\lfloor n/2\right\rfloor}/2^n$.
Therefore in this regime 
$\mathbb{P}(S_n = g) \leq \binom{n}{\left\lfloor n/2\right\rfloor}/2^n \sim
\frac{1}{\sqrt{n}}$. This shows that the inequality (\ref{TiepVU})
cannot be improved apart from the constant factor.
It is also very natural that the term $\frac{1}{m}$ is
dominant for $n \geq m^2$, exactly above the mixing time of $S_n$, that is known
to be of magnitude $m^2$ (see \cite{LevinPeresWilmer2006}, page 96).

In this paper we shall prove an optimal upper bound for $\rho(V_n)$, where the elements of the multiset $V_n$ lie in an arbitrary group. It turns out that a bound as in $(1)$ holds for arbitrary groups. Furthermore, for groups with elements having odd or infinite order we shall establish an optimal inequality for $\mathbb{P}(X_1\ast \cdots \ast X_n=x)$ 
without the requirement that the random variables $X_i$ are two-valued.

Let us remind the reader that we denote by $\varepsilon$ (usually supplied with a subscript) a uniform random variable on $\lbrace-1, 1\rbrace$. Sometimes it will be important to stress that these random variables are defined on $\mathbb{Z}_m$ instead of $\mathbb{R}$ and we shall do so on each occasion. 
We denote by $(a,b]_m$ and $[a,b]_m$ the set of integers in the intervals $(a,b]$ and $[a,b]$ modulo $m$. 
Given a natural number $m$, we shall write $\m$ for the smallest even
number such that $\m\geq m$. That is, we have $\m =
2\lceil\frac{m}{2}\rceil$.

\begin{theorem}\label{theo1}
Let $g_1, \ldots, g_n$ be elements of some group $G$ such that $|g_i|\geq m \geq 2$. Let $X_1,\ldots, X_n$ be independent random variables so that each $X_i$ has the uniform distribution on the two point set $\lbrace g^{-1}_i, g_{i} \rbrace$. Then for any $A\subset G$ with $|A|=k$ we have

\begin{equation}\label{ineq} \mathbb{P}\left(X_1\ast\cdots\ast X_n\in A\right)\leq \mathbb{P}\left(\varepsilon_1+\cdots+\varepsilon_n\in (-k,k]_{\m}\right),
\end{equation}

where $\varepsilon_i$ are independent uniform random variables on the set $\lbrace-1, 1\rbrace \subset \mathbb{Z}_{\m}$.\\

\end{theorem}

Note that Theorem~\ref{theo1} is optimal in the sense that if $G$ contains an element of order $\m$, the bound in $(2)$ can be attained. For instance, in the case $G=GL_{k}(\mathbb{C})$ the upper bound in $(2)$ is achieved by taking two point distributions concentrated on the diagonal matrix 
$\e^{\frac{2\pi i}{\m}}\mathbb{I}_{k}$ and its inverse. 
Theorem \ref{theo1} implies an inequality of the same type as the one by Tiep and Vu, but with a much better constant.
\begin{corollary}\label{corr}
Let $V_n=\lbrace g_1,\ldots, g_n \rbrace$ be  elements in some
group $G$ satisfying $|g_i|\geq m\geq 2$. Then 
\begin{equation}\label{corroll} 
\rho(V_n)\leq 
 \frac{2}{\m} + \sqrt{\frac{2}{\pi}}\frac{1}{\sqrt{n}}
 \leq
3\max\lbrace\frac{1}{m},\frac{1}{\sqrt{n}}\rbrace.
\end{equation}
\end{corollary}

The sequence of sums appearing on the right hand side of $(2)$ is a periodic Markov chain and so does not converge to a limit as $n\rightarrow \infty$. Nonetheless, it is well known that it does converge to a limit if we restrict the parity of $n$. Let us now express the quantity in the right hand side of $(2)$ in the case $|A|=1$ in asymptotic terms.

\begin{proposition}\label{prop}
Let $m\in \mathbb{N}$ and assume that $n\rightarrow \infty$. Then for any $l\in
\mathbb{Z}_{\m}$ of the same parity as $n$ we have
$$\mathbb{P}\left(\varepsilon_1+\cdots+\varepsilon_n = l\right)=\frac{2}{\m}+o(1).$$

\end{proposition}

The $o(1)$ term is actually exponentially small in terms of $n$. For such sharp 
quantitative estimates see \cite{Diaconis} pages 124-125. 
Note that Proposition~\ref{prop} implies that in (\ref{corroll}) the constant
after the last inequality cannot be smaller than $2$. Let us also note that both
constants in the expression 
 $\frac{2}{\m} + \sqrt{\frac{2}{\pi}}\frac{1}{\sqrt{n}}$
 are sharp. The term $\frac{2}{\m}$ is dominant in the case 
 $m,n \rightarrow \infty$ and $n \gg m^2$ and so Proposition~\ref{prop} shows
 that the constant $2$ cannot be reduced. In the case $m,n \rightarrow \infty$ 
and $n < m$ the therm $\sqrt{\frac{2}{\pi}}\frac{1}{\sqrt{n}}$ is dominating.
For $V_n = \{g, \dots, g\}$ for some element $g$ of order $\m$ we have
\begin{equation*}
  \rho(V_n) = \mathbb{P}\left(\varepsilon_1+\cdots+\varepsilon_n\in (-1,1]_{\m}\right) 
  = \frac{\binom{n}{\left\lfloor n/2\right\rfloor}}{2^n}
=(1+o(1))\sqrt{\frac{2}{\pi}}\frac{1}{\sqrt{n}}.
\end{equation*}

The simple random walk on $\mathbb{Z}_m$ for $m$ odd converges to the uniform
distribution on $\mathbb{Z}_m$ and so all probabilities converge to
$\frac{1}{m}$. It should now be unsurprising that the simple random walk on
$\mathbb{Z}_{m+1}$ is a much better "candidate" for a maximizer of the left hand
side in $(2)$, as by Proposition $\ref{prop}$ we gain an extra factor of $2$ asymptotically. 

From this point our prime focus will be on the particular case
$G=\mathbb{Z}^l_{m}$ for $m$ odd. In this case Theorem~\ref{theo1} does not provide the optimal bound. The approach we have for this case also works for certain groups other than $\mathbb{Z}^l_{m}$ and therefore we will state it in a general form. For $k \geq 1$ we define 
$$I_{n,k}^m=\left[\Bigl\lceil\frac{n-k+1}{2}\Bigl\rceil,\ldots,\Bigl\lceil\frac{n+k-1}{2}\Bigl\rceil\right]_m.$$
The latter set is an interval of $k$ points in $\mathbb{Z}_m$. We shall use the convention that $I_{n,0}^m=\emptyset$.

\begin{theorem}\label{theo2}
Let $X_1,\ldots, X_n$ be independent discrete random variables taking values in some group $G$ such that for each $i$ we have

\begin{equation}\label{half}
  \sup_{g\in G}\mathbb{P}\left(X_i=g\right)\leq \frac{1}{2}.
\end{equation}

Furthermore, assume that all non-identity elements in $G$ have odd or infinite order and that the minimal such order is at least some odd number $m\geq 3$. Then for any set $A\subset G$ of cardinality $k$ we have
\begin{equation*}
\mathbb{P}\left(X_1 \ast \cdots \ast X_n\in A\right) \leq  \mathbb{P}\left(\tau_1+\cdots+\tau_n\in I_{n,k}^m \right),
\end{equation*}
where $\tau_i$ are independent uniform random variables on the set $\lbrace 0, 1\rbrace \subset \mathbb{Z}_{m}$.
\end{theorem}

The distribution of $\tau_1+\cdots+\tau_n$ is asymptotically uniform in $\mathbb{Z}_m$ and thus we have $\mathbb{P}\left(X_1 \ast \cdots \ast X_n = g\right) \leq \frac{1}{m}+o(1)$. \\

\textbf{Remark 1.} Note that 
$$\mathbb{P}\left(\tau_1+\cdots+\tau_n\in I_{n,k}^m \right)=\mathbb{P}\left(\varepsilon_1+\cdots+\varepsilon_n\in 2I_{n,k}^m-n \right).$$
We formulated the result in terms of $\lbrace0,1\rbrace$ random variables $\tau_i$ for the sake of convenience only - in this formulation the set of maximum probability is an interval. As one notices, it is not so in formulating it in terms of $\lbrace-1,1\rbrace$ distributions $\varepsilon_i$.

\textbf{Remark 2.} The reason we restrict the elements to have odd order in
Theorem~\ref{theo2} is as follows. If there is an element of even order in the underlying group, then the group contains an element of order $2$, say $h$. 
Then by taking independent uniform random variables $X_i$ on the set $\lbrace 1,h \rbrace$ 
we obtain $\sup_{g\in G}\mathbb{P}\left(X_1 \ast \cdots \ast X_n= g\right)=\frac{1}{2}$.\\

In the case when $G$ is torsion-free we can actually prove that Erd\H{o}s's bound still holds even in this general setting.

\begin{proposition}
  Under the notation of Theorem~\ref{theo2} and assuming that $G$ is torsion-free for any set $A\subset G$ of cardinality $k$ we have

$$\mathbb{P}\left(X_1 \ast \cdots \ast X_n\in A\right) \leq  \mathbb{P}\left(\varepsilon_1+\cdots+\varepsilon_n\in (-k,k]\right),$$
where $\varepsilon_i$ are independent. In particular, for any $g\in G$ we have

$$\mathbb{P}\left(X_1 \ast \cdots \ast X_n = g\right)\leq \frac{\binom{n}{\left\lfloor n/2\right\rfloor}}{2^n}.$$

\end{proposition}
The latter proposition immediately follows by taking $m$ large enough in
Theorem~\ref{theo2} so that $\tau_1+\cdots+\tau_n$ is concentrated in a proper subset of $\mathbb{Z}_m$. For instance, assume that $m=n+2$. In this case the latter sum is strictly contained in $\mathbb{Z}_m$ and its probabilities are exactly the largest $k$ probabilities of $\varepsilon_1+\cdots+\varepsilon_n$ and we are done.

Our proofs are similar in spirit to Kleitman's approach in his solution of the
Littlewood-Offord problem in all dimensions. Actually, it is closer to a simplification of
Kleitman's proof in dimension 1 obtained in \cite{DJS}. The proofs thus proceed
by induction on dimension, taking into account a certain recurrence relation
satisfied by the worst-case random walk.  

\section{An open problem}

Theorem~\ref{theo1} gives an optimal inequality if an element with a given order
exists. To be more precise, if an element of order $\m$ exists. For groups
in which all elements have odd or infinite order, Theorem~\ref{theo2} gives the best possible result. It is thus natural to ask what happens if we have full knowledge of the orders of the elements of the underlying group $G$ and we are not in the aforementioned cases. The asymptotics of the cases when we do know the exact answer suggest the following guess.

\textbf{Conjecture}. Let $G$ be any group and fix an odd integer $m\geq 3$. Suppose that all possible even orders of elements in $G$ greater than $m$ are given by the sequence $S=\lbrace m_1, m_2, \ldots\rbrace$ in increasing order. Consider a collection of independent random variables $X_1,\ldots, X_n$ in $G$ such that each $X_i$ is concentrated on a two point set $\lbrace g_i, g^{-1}_i\rbrace$ and $|g_i|\geq m$. Then if $m_1<2m$ for any $A\subset G$ with $|A|=k$ we have
\begin{equation*} \mathbb{P}\left(X_1\ast\cdots\ast X_n\in A\right)\leq \mathbb{P}\left(\varepsilon_1+\cdots+\varepsilon_n\in (-k,k]_{m_1}\right),
\end{equation*}
where $\varepsilon_i$ are independent uniform random variables on the set $\lbrace-1, 1\rbrace \subset \mathbb{Z}_{m_1}$.\\ 
On the other hand, if $m_1\geq 2m$ we have   
\begin{equation*} \mathbb{P}\left(X_1\ast\cdots\ast X_n\in A\right)\leq \mathbb{P}\left(\tau_1+\cdots+\tau_n\in I_{n,k}^m\right),
\end{equation*}
where $\tau_i$ are independent uniform random variables on the set $\lbrace 0, 1\rbrace \subset \mathbb{Z}_{m}$.\\

If true, the latter conjecture would settle the remaining cases.

\section{Proofs}

In order to prove Theorems~\ref{theo1}-\ref{theo2}, we shall require a simple group theoretic statement contained in the following lemma.

\begin{lemma}\label{lemma1}
Let $G$ be a group and $g\in G$ be an element of order greater then or equal to $m\geq 2$. Then for any
finite set $A\subset G$ and a positive integer $s$ 
such that $s < \frac{m}{|A|}$ we have $A\neq Ag^{s}$. 
\end{lemma}

\textit{\textbf{Proof of Lemma~\ref{lemma1}.}}
Suppose there is a nonempty set $A \subset G$ and a positive integer $s$ such
that $|A| = k < \frac{m}{s}$ and $A = Ag^s$. Take some $a \in A$ and consider
elements $ag^{si}$, $i=0 \dots k$. All these $k+1$ elements are in
the set $A$ hence at least two of them must be equal. Let us say $ag^{s i} = ag^{s j}$
for some integers $0 \leq i < j \leq k$. But this immediately gives
a contradiction since then $g^{s  (j-i)}$ is equal
to the group identity element and  ${m \leq s (j-i) \leq s k}$.

\textit{\textbf{Proof of Theorem~\ref{theo1}.}}
If $n=1$ the inequality (\ref{ineq}) is trivial. For $k\geq \frac{m}{2}$ and all $n$ the right hand side of (\ref{ineq}) becomes $1$ since in this case $(-k,k]_{\m}$ covers the support of the sum $\varepsilon_1+\cdots+\varepsilon_n$ and so there is nothing to prove. We shall henceforth assume that $n>1$ and $k< \frac{m}{2}$.\\
By Lemma~\ref{lemma1} we have that $Ag_{n}\neq Ag^{-1}_{n}$. Take some $h\in
Ag_{n} \backslash Ag^{-1}_{n}$ and define ${B=Ag_{n}\backslash \lbrace h \rbrace}$ and $C=Ag^{-1}_{n} \cup \lbrace h \rbrace$. We then have

\begin{eqnarray}
&&2 \mathbb{P}(X_1 \ast \cdots \ast X_n\in A)=\mathbb{P}(X_1 \ast \cdots \ast X_{n-1}\in Ag_n)+\mathbb{P}(X_1 \ast \cdots \ast X_{n-1}\in Ag^{-1}_n) \nonumber\\
&=&\mathbb{P}(X_1 \ast \cdots \ast X_{n-1}\in B)+\mathbb{P}(X_1 \ast \cdots \ast X_{n-1}\in C)\\
&\leq &\mathbb{P}(\varepsilon_1+\cdots+\varepsilon_{n-1}\in (-k-1,k+1]_{\m})+\mathbb{P}(\varepsilon_1+\cdots+\varepsilon_{n-1}\in (-k+1,k-1]_{\m})\\
&=&\mathbb{P}(\varepsilon_1+\cdots+\varepsilon_{n-1}\in (-k-1,k-1]_{\m})+\mathbb{P}(\varepsilon_1+\cdots+\varepsilon_{n-1}\in (-k+1,k+1]_{\m})\\
&=&2\mathbb{P}\left(\varepsilon_1+\cdots+\varepsilon_n\in (-k,k]_{\m}\right)\nonumber.
\end{eqnarray}
This completes the proof.

\textbf{Remark 3}. Note that in $(6)$-$(7)$ we used the fact that for $k<\frac{m}{2}$ the sets $(-k+1,k-1]_{\m}$ and $(k-1,k+1]_{\m}$ are disjoint in $\Z_{\m}$.

In the proof of Theorem~\ref{theo2} we shall make use of the following simple
lemma which will allow us to switch from general distributions satisfying the
condition (\ref{half}) to two-point distributions.

\begin{lemma}\label{lemma2}
Let $X$ be a random variable on some group $G$ that takes only finitely many values, say $x_1,\ldots, x_n$. Suppose that $p_i=\mathbb{P}(X=x_i)$ are rational numbers and that $p_i\leq \frac{1}{2}$. Then we can express the distribution of $X$ as a convex combination of distributions that are uniform on some two point set.

\end{lemma}

\textit{\textbf{Proof of Lemma~\ref{lemma2}.}} Denote by $\mu$ the distribution of $X$. Since the $p_i$'s are all rational, we have $p_i=\frac{k_i}{K_i}$ for some $k_i,K_i\in \mathbb{Z}$. We shall now view $\mu$ as a distribution on a multiset $M$ made from the elements $x_i$ in the following way - take $x_i$ exactly $2 k_i\prod_{j\neq i}K_i$ times into $M$. This way $\mu$ has the uniform distribution on $M$. We thus have that $M=\lbrace y_1, \ldots, y_{2N}\rbrace$ for the appropriate $N$. Construct a graph on the elements on $M$ by joining two of them by an edge if and only if they are distinct. Since we had $p_i\leq \frac{1}{2}$, each vertex of this graph has degree at least $N$. Thus by Dirac's Theorem, our graph contains a Hamiltonian cycle, and, consequently - a perfect matching. Let $\mu_i$ be the uniform distribution on two vertices of the latter matching ($i=1,2,\ldots, N$). We have 

$$\mu=\frac{1}{N}\sum_{i=1}^{N}\mu_{i}.$$

\textit{\textbf{Proof of Theorem~\ref{theo2}.}} We shall argue by induction.
First notice that the claim of the Theorem is true for $n=1$. Furthermore, it is
also true for $k\geq m$ since in that case the bound for the probability in
question becomes $1$. We therefore shall from now on assume that $n>1$ and
$1\leq k \leq m-1$.  Denote by $\mu_i$ the distribution of the random variable
$X_i$. We can without loss of generality assume that each $X_i$ is concentrated
on finitely many points and that for each $g\in G$ we have $\mathbb{P}(X_i =
g)\in \mathbb{Q}$. By Lemma~\ref{lemma2}, each $\mu_i$ can be written as a
convex combination of distributions that are uniform on some two-point set.
Define the random variable $f_{i}(X_i)= \mathbb{E}_i1\lbrace X_1 \ast \cdots \ast X_n \in A \rbrace$, where $\mathbb{E}_i$ stands for integration with respect to all underlying random variables except $X_i$. Then for each $i$ we have
\begin{equation}
\mathbb{P}\left(X_1 \ast \cdots \ast X_n\in A\right)=\mathbb{E}f_{i}(X_i).
\end{equation}
The latter expectation is linear with respect to the distribution of $X_i$. Therefore we can assume that it will be maximized by some choice of two-point distributions coming from the decomposition of $\mu_i$. We shall therefore from this point assume that $X_n$ takes only two values, say $h_1$ and $h_2$, with equal probabilities.

Note that the intervals $I_{n,k}^m$ have recursive structure. Namely, if $1\leq k \leq m-1$ and we regard them as multisets, we have the relation $I_{n,k}^m\cup (I_{n,k}^m-1)= I_{n-1,k-1}^m \cup I_{n-1,k+1}^m$. The pairs on intervals appearing on both sides of the latter equality heavily overlap. This means that we can take one endpoint of $I_{n-1,k+1}^m$ that does not belong to $I_{n-1,k-1}^m$ and move it to this shorter interval. The resulting intervals are both of length $k$ and are exactly the intervals $I_{n,k}^m$ and $I_{n,k}^m-1$. We shall use this after the inductive step. 

Take a finite set $A\subset G$ with $k$ elements. Note that the element
$h^{-1}_{2}h_{1}\neq 1_G$ and so it has order at least $m$. By
Lemma~\ref{lemma1} we have that $Ah^{-1}_1\neq Ah^{-1}_2$ as $A\neq Ah^{-1}_{2}h_{1}$. Take some $h\in Ah^{-1}_1 \backslash Ah^{-1}_2$ and define $B=Ah^{-1}_1 \backslash \lbrace h \rbrace$ and $C=Ah^{-1}_2 \cup \lbrace h \rbrace$. We have

\begin{eqnarray*}
2 \mathbb{P}(X_1 \ast \cdots \ast X_n\in A)&=&\mathbb{P}(X_1 \ast \cdots \ast X_{n-1}\in Ah^{-1}_1)+\mathbb{P}(X_1 \ast \cdots \ast X_{n-1}\in Ah^{-1}_2)\\
&=&\mathbb{P}(X_1 \ast \cdots \ast X_{n-1}\in B)+\mathbb{P}(X_1 \ast \cdots \ast X_{n-1}\in C)\\
&\leq & \mathbb{P}(\tau_1 + \cdots + \tau_{n-1}\in I_{n-1,k-1}^m)+\mathbb{P}(\tau_1 + \cdots + \tau_{n-1}\in I_{n-1,k+1}^m)\\
&=& \mathbb{P}(\tau_1 + \cdots + \tau_{n-1}\in I_{n,k}^m-1)+\mathbb{P}(\tau_1 + \cdots + \tau_{n-1}\in I_{n,k}^m)\\
&=&2\mathbb{P}(\tau_1 + \cdots + \tau_{n}\in I_{n,k}^m).
\end{eqnarray*}
This completes the proof.
\\

\textit{\textbf{Proof of Corollary~\ref{corr}.}} 
We shall use an identity on evenly spaced binomial coefficients proved in \cite{BinomSum}:
\begin{equation}\label{binomeq} 
  \binom{n}{t} + \binom{n}{t+s} + \binom{n}{t+2s} + \dots =
  \frac{1}{s}\sum_{j=0}^{s-1}{\Big(2\cos{\frac{i\pi}{s}}\Big)^n\cos{\frac{\pi(n-2t)j}{s}}}.
\end{equation}

By Theorem~$\ref{theo1}$ we have 
\begin{equation}\label{coreq} 
\rho(V_n) \leq 
\mathbb{P}\left(\varepsilon_1+\cdots+\varepsilon_n\in (-1,1]_{\m}\right).
\end{equation}

The right hand of the equation (\ref{coreq}) is the sum of binomial
probabilities $\binom{n}{i}/2^{n}$,
where $i$ is such that $2i-n$ is congruent to $1_{\{n \in 2\Z + 1\}}$ modulo $\m$.
Let $t$ be the residue of ${(n- 1_{\{n \in 2\Z + 1\}})/2}$
modulo $\frac{\m}{2}$. 

Using the identity (\ref{binomeq}) and the elementary inequalities 
 $\cos x\leq \exp(-x^2/2)$ for $x \in [0, \frac{\pi}{2}]$ and
$\int_{0}^{\infty}\e^{\frac{-x^2}{2\sigma^2}}dx \leq \frac{\sigma\sqrt{2\pi}}{2}$
we obtain

\begin{eqnarray}
\mathbb{P}\left(\varepsilon_1+\cdots+\varepsilon_n\in
(-1,1]_{\m}\right) &=&
\frac{
  \binom{n}{t} + \binom{n}{t+\m/2} + \binom{n}{t+2\m/2} + \dots
}{2^n} \nonumber \\
&=&
\frac{2}{\m}\sum_{j=0}^{\frac{\m}{2}-1}{\Big(2\cos{\frac{2i\pi}{\m}}\Big)^n
\cos{\frac{2\pi(n-2t)j}{\m}}} \nonumber \\
&\leq&
\label{costrick}
\frac{2}{\m} + \frac{2}{\m}\sum_{j=1}^{\frac{\m}{2}-1
}\Big|\cos{\frac{2j\pi}{\m}}\Big|^n \\
&\leq&
\frac{2}{\m} + \frac{4}{\m}\sum_{j=1}^{\left\lfloor \frac{\m}{4} \right\rfloor
}\Big|\cos{\frac{2j\pi}{\m}}\Big|^n \nonumber \\
&\leq& 
\frac{2}{\m} + \frac{4}{\m}\sum_{j=1}^{\left\lfloor \frac{\m}{4} \right\rfloor}
\e^{-2\pi^2j^2n/\m^2} \nonumber \\
&<&\frac{2}{\m} + \frac{4}{\m} \int_{0}^{\infty}\e^{-2\pi^2x^2n/\m^2} dx
\nonumber \\
&\leq& \frac{2}{\m} + \sqrt{\frac{2}{\pi}}\frac{1}{\sqrt{n}} 
\leq \frac{2}{m} + \sqrt{\frac{2}{\pi}}\frac{1}{\sqrt{n}} \nonumber .
\end{eqnarray}

Note that in (\ref{costrick}) we replaced $|\cos{\frac{2\pi j}{\m}}|$ by $|\cos{\frac{\pi(\m - 2j)}{\m}}|$ when
$j > \frac{\m}{4}$.
This completes the proof.

\bibliographystyle{plain}

\end{document}